\documentclass[fleqn,11pt]{article}

\usepackage{amsmath, amscd, amsfonts, amssymb, graphicx, color}
\usepackage{hyperref}
\usepackage{cite}
\usepackage{amsmath}
\usepackage{amsfonts}
\usepackage{epsfig}
\usepackage{amssymb}
\usepackage{amsthm}
\usepackage{epstopdf}
\usepackage{footnote}
\usepackage{newlfont}
\usepackage{color}
\newtheorem{theorem}{Theorem}[section]
\newtheorem{lemma}[theorem]{Lemma}

\newtheorem{corollary}[theorem]{Corollary}
\theoremstyle{definition}
\newtheorem{definition}[theorem]{Definition}
\theoremstyle{example}
\newtheorem{example}[theorem]{Example}

\newtheorem{remark}[theorem]{Remark}

\begin{document}

\title{On the distance from a matrix polynomial to matrix polynomials with two prescribed eigenvalues}


 \author{E. Kokabifar\thanks{Department of Mathematics,
 Faculty of Science, Yazd University, Yazd, Iran
 (e.kokabifar@stu.yazd.ac.ir, loghmani@yazd.ac.ir).},\,
  G.B. Loghmani,\footnotemark[1]\,\,
   A.M. Nazari\thanks{Department of Mathematics, Faculty of Science, Arak University, Arak, Iran (a-nazari@araku.ac.ir).},\, and S.M. Karbassi\thanks{Department of Mathematics, Yazd Branch, Islamic Azad University, Yazd, Iran (mehdikarbassi@gmail.com).}}

\maketitle

\vspace{-6mm}

\begin{abstract}
Consider an $n \times n$ matrix polynomial $P(\lambda)$. A spectral norm distance from $P(\lambda)$ to the set of
$n \times n$ matrix polynomials that have a given scalar
$\mu\in\mathbb{C}$ as a multiple eigenvalue was introduced and obtained by Papathanasiou and Psarrakos. They computed lower and upper bounds for this distance, constructing an associated perturbation of $P(\lambda)$. In this paper, we extend this result to the case of two given distinct complex numbers $\mu_{1}$ and $\mu_{2}$. First, we compute a lower bound for the spectral norm distance from $P(\lambda)$ to the set of matrix polynomials that have $\mu_1,\mu_2$ as two eigenvalues. Then we construct an associated perturbation of $P(\lambda)$, such that the perturbed matrix polynomial has two given scalars $\mu_1$ and $\mu_2$ in its spectrum. Finally, we derive an upper bound for the distance by the constructed perturbation of $P(\lambda)$. Numerical examples are provided to illustrate the validity of the method.
\end{abstract}



\section{Introduction and preliminaries}

\noindent

Let $A$ be an $n\times n$ complex matrix and  let $L$ be the set of complex $n\times n$ matrices with a multiple zero eigenvalue. In 1999, Malyshev \cite{malyshev} obtained a formula for the spectral norm distance from $A$ to $L$ which can be considered as a theoretical solution to
Wilkinson's problem, that is, the calculation of the distance from a matrix $A \in \mathbb{C}^{n \times n}$ that has all its eigenvalues simple to the $n\times n$ matrices with multiple
eigenvalues. Wilkinson introduced this distance in
\cite{wilkinson}, and some bounds for it were computed by Ruhe \cite{ruhe}, Wilkinson \cite{wil1,wil2,wil3,wil4} and Demmel
\cite{demmel1}. Also, Malyshev's results were extended by Lippert \cite{lipert} and Gracia \cite{gracia}; they obtained a spectral norm distance from $A$ to the set of matrices that have two prescribed eigenvalues and studied a nearest matrix with the two desired eigenvalues. 
In 2008, Papathanasiou and Psarrakos \cite{papa} introduced and
studied a spectral norm distance from a $n\times n$ matrix
polynomial $P(\lambda)$ to the set of $n\times n$ matrix
polynomials that have a scalar $\mu\in\mathbb{C}$ as a multiple
eigenvalue. In particular, generalizing Malyshev's methodology,
they computed lower and upper bounds for this distance,
constructing an associated perturbation of $P(\lambda)$ for the
upper bound.

In this paper, motivated by the above, extending some of the results obtained in\cite{papa} for the case of two distinct eigenvalues is considered. This note concerns the bounds for a spectral norm distance from an $n \times n$ matrix polynomial $P(\lambda)$ to the set $\mathcal{P}_{\mu_1,\mu_2}$ of $n\times n$ matrix polynomials that have two given distinct scalars $\mu_1, \mu_2\in\mathbb{C}$ in their spectrum. In addition, construction of an associated perturbation of $P(\lambda)$ is also considered. Replacing the divided differences by derivative of $P(\lambda)$ in \cite[Definition 5]{papa}, extending all of necessary definitions and lemmas in \cite{papa, gracia, lipert, malyshev}, and also constructing an appropriate perturbation of $P(\lambda)$ are some of the main ideas used in this article. This paper can be considered as generalization of the results obtained in \cite{lipert, gracia} for the case of matrix polynomials. In the next section, some definitions for a matrix polynomial presented and also a spectral norm distance from $P(\lambda)$ to $\mathcal{P}_{\mu_1,\mu_2}$ is introduced. In Section 3, we prove some lemmas which will be applied in the Section 4 where we derive lower and upper bounds of $P(\lambda)$ for the distance from $P(\lambda)$ to $\mathcal{P}_{\mu_1,\mu_2}$ and construct an associated perturbation of $P(\lambda)$. In Section 5, connection between the previous result and ours is discussed. Finally, in last section, a numerical example is given to illustrate the validation and application of our method.

Consider an $n \times n$ matrix polynomial
\begin{equation}\label{plambda}
P(\lambda ) = A_m \lambda ^m  + A_{m - 1} \lambda ^{m - 1}  + ...
+ A_1 \lambda  + A_0,
\end{equation}
where $A_j\in\mathbb{C}^{n \times n}(j = 0,1,...,m)$ with det$(A_m)\neq 0$ and $\lambda$ is a complex variable.
The study of matrix polynomials, especially with regard to their
spectral analysis, has received a great deal of attention and has
been used in many applications \cite{glr,kacz,lanc,markus}.
Standard references for the theory of matrix polynomials are
\cite{glr,markus}. Here, some definitions of matrix polynomials
are briefly reviewed.

If for a scalar $\lambda_0 \in \mathbb{C}$ and some nonzero vector
$x_0 \in {\mathbb{C}^{n}}$, it holds that $P(\lambda_0) x_0 = 0$,
then the scalar $\lambda_0$ is called an \textit{eigenvalue} of
$P(\lambda)$ and the vector $x_0$ is known as a \textit{(right)
eigenvector} of $P(\lambda)$ corresponding to $\lambda_0$. The
\textit{spectrum} of $P(\lambda)$, denoted by $\sigma(P)$, is the
set of all eigenvalues of $P(\lambda)$. Since the leading
matrix-coefficient $A_m$ is nonsingular, the spectrum $\sigma(P)$
contains at most $mn$ distinct finite elements. The multiplicity
of an eigenvalue $\lambda_0 \in \sigma(P)$ as a root of the scalar
polynomial $\det P(\lambda)$ is said to be the \textit{algebraic
multiplicity} of $\lambda_0$, and the dimension of the null space
of the (constant) matrix $P(\lambda_0)$ is known as the
\textit{geometric multiplicity} of $\lambda_0$. Algebraic
multiplicity of an eigenvalue is always greater than or equal to
its geometric multiplicity. An eigenvalue is called
\textit{semisimple} if its algebraic and geometric multiplicities
are equal, otherwise it is known as \textit{defective}.
\begin{definition}
Let $P(\lambda )$ be a matrix polynomial as in (1) and let $\Delta _j  \in \mathbb{C}^{n \times n}, (j =
0,1,...,m)$ be the arbitrary matrices. We consider perturbations of the matrix polynomial $P(\lambda)$ as following
\begin{equation}
Q(\lambda ) = P(\lambda ) + \Delta (\lambda ) = \sum\limits_{j =
0}^m {(A_j  + \Delta _j )\lambda ^j }.
\end{equation}
Moreover, consider the positive quantity $\varepsilon>0$ and set of given weights $w = \{ \omega _0 ,\omega _1 ,...,\omega _m \}$, such that $w$ is a set of nonnegative coefficients with $\omega _0 > 0$. Define the associated set of perturbations of $P(\lambda )$ by
\begin{equation*}
{\mathcal{B}}(P,\varepsilon ,w) = \{ Q(\lambda ){  ~ as ~in ~ }
(2): \left\| {\Delta _j } \right\| \le \varepsilon \omega
_j,~~~ j = 0,1,...,m\},
\end{equation*}
and consider the scalar polynomial $w(\lambda )$ corresponding to the weights as following
\begin{equation*}
w(\lambda ) = \omega _0+ \omega _1 \lambda + \ldots  + \omega _{m - 1} \lambda ^{m- 1}  + \omega _m \lambda ^m.
\end{equation*}
\end{definition}
\begin{definition}\label{dis}
Let the matrix polynomial $P(\lambda )$  as in
(\ref{plambda}) and two distinct complex numbers $\mu _1$ and $\mu_2$ are given. Define the distance from $P(\lambda)$ to $\mathcal{P}_{\mu_1,\mu_2}$ by
\begin{equation*}
\mathcal{D} (P,\mu _1 ,\mu _2 ) = \min \{ \varepsilon  \ge
0~:~\exists Q(\lambda ) \in {\mathcal{B}}(P,\varepsilon ,w) { ~with
~\mu_1 ~and~\mu_2 ~as~ two~ eigenvalues}  \}.
\end{equation*}
\end{definition}
\begin{definition}
Let $P(\lambda )$ be a matrix polynomial as in (\ref{plambda}) and $\mu_1$ and $\mu_2$ be two given
distinct complex numbers. Define the $2n \times 2n$ matrix
\begin{equation*}
F[P(\mu _1 ,\mu _2 );\gamma ] = \left(\left[
{\begin{array}{*{20}c}
   {P(\mu _1 )} & 0  \\
   {\gamma \frac{{P(\mu _1 ) - P(\mu _2 )}}{{\mu _1  - \mu _2 }}} & {P(\mu _2 )}  \\
\end{array}}\right] \right);\qquad \gamma  \in \mathbb{C}.
\end{equation*}
\end{definition}
Henceforth for simplicity we denote $\frac{{P(\mu _1
) - P(\mu _2 )}}{{\mu _1  - \mu _2 }}$ by $P[\mu _1 ,\mu _2
]$ and so on.
\section{Properties of $s_{2n-1}(F[P(\mu _1 ,\mu _2 );\gamma ])$ and its corresponding singular vectors}
In this section we study some properties of $s_{2n-1}(F[P(\mu _1 ,\mu _2 );\gamma ])$ and its corresponding singular vectors. These properties are needed in the next section in order to obtain bounds for $\mathcal{D} (P,\mu _1 ,\mu _2 )$ and construct a perturbation of $P(\lambda)$. In this section some definitions and lemmas of [4-6] are reconstructed for the case of two distinct eigenvalues. Proving of some lemmas is mostly similar to the proof of related lemmas in its references. Therefore, for convince, this proofs can be omitted.
\begin{lemma}
For $\mu_1, \mu_2 \in \mathbb{C}$ and for all $\gamma \neq 0 $, we have either $s_{2n -1} (F[P(\mu _1 ,\mu _2 ); \gamma ]) \neq 0 $ or $s_{2n - 1} (F[P(\mu _1 ,\mu _2 );\gamma ]) \equiv 0$ .\\
\end{lemma}
\textbf{Proof.} Similar to Lemma 3.4 of \cite {lipert} can be
verified easily.\qquad $\square$
\begin{lemma}\label{lem2}
If $\mu_1$ and $\mu_2$ are two eigenvalues of the
matrix polynomial $ Q(\lambda ) =
P(\lambda ) + \Delta (\lambda ) $, then for any $\gamma \neq 0$
\begin{equation*}
s_{2n - 1} (F[P(\mu _1 ,\mu _2 );\gamma ]) \le \left\| {F[\Delta
(\mu _1 ,\mu _2 );\gamma ])} \right\|.
\end{equation*}
\end{lemma}
\textbf{Proof.} Let $\mu_1$ and $\mu_2$ be two eigenvalues of $ Q(\lambda ) =P(\lambda ) + \Delta (\lambda )
$, then for any $\gamma \neq 0$
\begin{equation*}
s_{2n - 1} (F[Q(\mu _1 ,\mu _2 );\gamma ]) = s_{2n} (F[Q(\mu _1
,\mu _2 );\gamma ]) = 0,
\end{equation*}
applying the Weyl inequalities for singular values (for example, see Corollary 5.1 of \cite{demmel})  for the above relation yields
\begin{equation*}
\left| {s_{2n - 1} (F[Q(\mu _1 ,\mu _2 );\gamma ]) - s_{2n - 1}
(F[P(\mu _1 ,\mu _2 );\gamma ])}
 \right| \le \left\| {F[\Delta (\mu _1 ,\mu _2 );\gamma ]}
 \right\|,
\end{equation*}
combining two recent relation concludes
\begin{equation*}
s_{2n - 1} (F[P(\mu _1 ,\mu _2 );\gamma ])\le \left\| {F[\Delta
(\mu _1 ,\mu _2 );\gamma ]} \right\|.\qquad \square
\end{equation*}
The two above lemmas will be used to obtain a lower bound for $\mathcal{D} (P,\mu _1 ,\mu _2 )$. In remainder of this section, some properties of singular vectors of $s_{2n - 1} (F[P(\mu _1 ,\mu _2 );\gamma ])$ are studied which will be necessary for computation an upper bound for $\mathcal{D} (P,\mu _1 ,\mu _2 )$ and a perturbation of $P(\lambda)$ in next section.
\begin{definition}\label{UV}
Let $\left[ {\begin{array}{*{20}c}
  {u_1 (\gamma )}  \\
   {u_2 (\gamma )}  \\
\end{array}} \right],\left[ {\begin{array}{*{20}c}
   {v_1 (\gamma )}  \\
  {v_2 (\gamma )}  \\
\end{array}} \right] \in \mathbb{C}^{2n} (u_k (\gamma ),v_k (\gamma ) \in \mathbb{C}^n,k = 1,2)
$ be a pair of left and right singular vectors of $ s_{2n - 1}
(F[P(\mu _1 ,\mu _2 );\gamma ])$ respectively. Define the $n \times 2$ matrices $U(\gamma )
= [u_1 (\gamma )~u_2(\gamma )],$ and $ V(\gamma ) = [v_1 (\gamma
)~v_2 (\gamma )]$.

Now set $\hat{u} (\gamma ) = u_2 (\gamma ) - \theta u_1 (\gamma),  \hat{v} (\gamma ) = v_2 (\gamma ) - \theta v_1 (\gamma)$ where $\theta  = \frac{\gamma }{{\mu _1  - \mu _2 }}$ and define $\hat{U}(\gamma )
= [u_1 (\gamma )~\hat{u}(\gamma )],$ and $ \hat{V}(\gamma )
= [v_1 (\gamma )~\hat{v}(\gamma )]$.
\end{definition}
To construct a perturbation of $P(\lambda)$ that has $\mu_1$ and $\mu_2$ as two eigenvalues and obtain a upper bound, we need have rank$(\hat{V}(\gamma ))=2.$ Now we derive a sufficient condition that implies it.  It is easy to verify that $ s_{2n - 1}
(F[P(\mu _1 ,\mu _2 );\gamma ])$ is an even function of $\gamma$, therefore, without loss of generality, hereafter we can assume that the parameter $\gamma$ is a nonnegative real number. Note that the following lemma which can be verified easily by considering Lemma 3 of \cite{gracia} and Lemma 13 of \cite{papa} provides a condition that assures that the function  $s_{2n - 1} (F[P(\mu _1 ,\mu_2 );\gamma ])$ attains its maximum value at a finite point.
\begin{lemma}
If  rank$(P[\mu _1 ,\mu _2
]) \ge 2$, Then $\mathop {\lim }\limits_{\gamma  \to \infty } s_{2n - 1} (F[P(\mu _1 ,\mu
_2 );\gamma ]) = 0.$
\end{lemma}

This corollary concludes that if rank$(P[\mu _1 ,\mu _2]) \ge 2$, then there is a finite point $\gamma_*\geq0$ where the singular value $s_{2n - 1} (F[P(\mu _1 ,\mu_2 );\gamma ])$ attains its maximum. Henceforth, for the sake of simplicity, $s_*$ denotes this maximum value of $s_{2n - 1} (F[P(\mu _1 ,\mu_2 );\gamma ])$, i.e.,
\begin{equation*}
s_*  = \mathop {\max }\limits_{\gamma  \ge 0} s_{2n - 1} (F[P(\mu
_1 ,\mu _2 );\gamma ]) = s_{2n - 1} (F[P(\mu _1 ,\mu _2 );\gamma
_* ]),
\end{equation*}
and $\theta_*  = \frac{\gamma_* }{{\mu _1  - \mu _2 }}$. It is obvious that if $s_* =0$, then $\mu_1$ and $\mu_2$ are two eigenvalues of $P(\lambda)$. Therefore, in what follows we assume that $s_*>0$.

By applying the lemma 5 of \cite{malyshev} for $F[P(\mu _1 ,\mu _2 );\gamma]$ we
have the next result.
\begin{lemma}\label{lemma6}
Let $\mu_1$ and $\mu_2$ be two complex numbers and let
$\gamma_*>0$. Then there exist a pair $\left[
{\begin{array}{*{20}c}
   {u_1 (\gamma _* )}  \\
   {u_2 (\gamma _* )}  \\
\end{array}} \right], \left[ {\begin{array}{*{20}c}
   {v_1 (\gamma _* )}  \\
   {v_2 (\gamma _* )}  \\
\end{array}} \right] \in\mathbb{C} ^{2n}~( u_k (\gamma _* ),v_k (\gamma _* ) \in \mathbb{C}^n,~k =
1,2)$ of left and right singular vectors of $s_*$ respectively,
such that

$1.~u_2 (\gamma _* )^*P[\mu _1 ,\mu _2 ]v_1 (\gamma _* ) =0, $

$2.~u_2 (\gamma _* )^*u_1 (\gamma _* ) = v_2 (\gamma _*)^*v_1(\gamma _* )$, and

3.~for the $n\times 2$ matrices $ U(\gamma _* ) = [u_1 (\gamma _*
)~u_2 (\gamma _* )]_{n \times 2}$ and $ V(\gamma _* ) = [v_1
(\gamma _* )~v_2 (\gamma _* )]_{n \times 2} $ we have $
U(\gamma _* )^* U(\gamma _* ) = V(\gamma _* )^* V(\gamma _* ).$
\end{lemma}
\textbf{Proof.} The first part of proof is similar to the first part of the proof of \cite[Lemma 17]{papa}. For the second part, we know that the vectors $\left[
{\begin{array}{*{20}c}
   {u_1 (\gamma _* )}  \\
   {u_2 (\gamma _* )}  \\
\end{array}} \right], \left[ {\begin{array}{*{20}c}
   {v_1 (\gamma _* )}  \\
   {v_2 (\gamma _* )}  \\
\end{array}} \right]$ satisfy the following relations
\begin{equation}\label{fn}
F[P(\mu _1 ,\mu _2 );\gamma _* ]\left[ {\begin{array}{*{20}c}
   {v_1 (\gamma _* )}  \\
   {v_2 (\gamma _* )}  \\
\end{array}} \right] = s_* \left[ {\begin{array}{*{20}c}
   {u_1 (\gamma _* )}  \\
   {u_2 (\gamma _* )}  \\
\end{array}} \right],
\end{equation}\vspace{-.1cm}
\begin{equation}\label{fs}
\left[ {\begin{array}{*{20}c}
   {u_1 (\gamma _* )^* }  \\
   {u_2 (\gamma _* )^* }  \\
\end{array}} \right]F[P(\mu _1 ,\mu _2 );\gamma _* ] = s_* \left[ {\begin{array}{*{20}c}
   {v_1 (\gamma _* )^* }  \\
   {v_2 (\gamma _* )^* }  \\
\end{array}} \right],
\end{equation}
now from equations (\ref{fn}), (\ref{fs}) and $u^* _2 (\gamma _* )P[\mu _1 ,\mu _2 ]v_1 (\gamma _* ) = 0,$ we have 
\[
\begin{array}{*{20}c}
   {u_2 (\gamma _* )^* \left( {P(\mu _1 ) - P(\mu _2 )} \right)v_1 (\gamma _* ) = 0, \Leftrightarrow }  \\
   {u_2 (\gamma _* )^* \left( {P(\mu _1 )v_1 (\gamma _* )} \right) - \left( {u_2 (\gamma _* )^* P(\mu _2 )} \right)v_1 (\gamma _* ) = 0, \Leftrightarrow }  \\
   {u_2 (\gamma _* )^* \left( {s_* u_1 (\gamma _* )} \right) - \left( {s_* v_2 (\gamma _* )^* } \right)v_1 (\gamma _* ) = 0, \Leftrightarrow }  \\
   {s_* \left( {u_2 (\gamma _* )^* u_1 (\gamma _* ) - v_2 (\gamma _* )^* v_1 (\gamma _* )} \right) = 0.}  \\
\end{array}
\]
Since $s_*>0$, thus $u_2 (\gamma _* )^*u_1 (\gamma _* ) = v _2 (\gamma _* )^*v_1(\gamma _* ).$ Third part of proof can be verified be considering the two previous parts of proof and following the procedure described for the second part of proof of \cite[Lemma 17]{papa}.

Note that third part of the above lemma deduces that if there exist a linear combination  of $u_1(\gamma_*)$ and $u_2(\gamma_*)$, then simultaneously we have the same linear combination of $v_1(\gamma_*)$ and $v_2(\gamma_*)$.
\begin{corollary}\label{natijeh}
Two matrices $ \hat{U}(\gamma
_* )$ and $ \hat{V}(\gamma_* )$, satisfy $\hat{U} (\gamma _* )^*\hat{U}(\gamma _* ) = \hat{V} (\gamma _*)^*\hat{V}(\gamma _* ).$
\end{corollary}
The following lemma provides a sufficient condition implying rank$(\hat V(\gamma _* )) = 2$.
\begin{lemma}\label{lemma7}
Suppose that $\gamma_*>0$ and $P[\mu _1 ,\mu _2 ]$ is a nonsingular
matrix. Then the two matrices $U({\gamma _*})$ and $V({\gamma _*})$ are full rank.
\end{lemma}
\textbf{Proof.} First it is shown that the four vectors $u_1 (\gamma _* ), u_2 (\gamma _* ), v_1(\gamma _* )~ and~v_2 (\gamma _* )$ are nonzero vectors. Assume the contrary, for example, that $u_2(\gamma_*)=0.$ Then the third part of Lemma \ref{lemma6} implies $v_2(\gamma_*)=0,$ and also the equation (\ref{fn}) yields $\gamma_*P[\mu _1 ,\mu _2 ]v_1(\gamma_*)=0$. Since $\gamma_*>0$ and $P[\mu _1 ,\mu _2 ]$ is an invertible matrix, we derive $v_1(\gamma_*)=0$ which is a contradiction because $v_1(\gamma_*), v_2(\gamma_*)$ form the right singular vector of $s_*$. By following the same reasoning that is used for $u_2(\gamma_*)$, we can derive that remainder of vectors are also nonzero.

Now, we will prove that $\hat{U}(\gamma_*)$ is a full rank matrix. Clearly, this concludes that $\hat{u} (\gamma )$ is a nonzero vector. Since $P[\mu _1 ,\mu _2 ]$ is nonsingular, then $M=\rho \left(P(\mu_1)-P(\mu_2)\right)$ for any $\rho \ne 0,$ is also a nonsingular matrix. Suppose from the contrary that for a nonzero $\xi\in \mathbb{C}$ we have $u_2 (\gamma _* )=\xi u_1 (\gamma _*)$. Two cases are considered.

\textit{Case 1.}\ Consider the case for which $\xi\ne\theta_*.$ Then from (\ref{fn}) we obtain
\begin{equation}\label{n1}
{P(\mu _1 )v_1 (\gamma _* ) = s_* u_1 (\gamma _* )},
\end{equation}
and
\begin{equation}\label{n2}
{\theta_* P(\mu _1 )v_1 (\gamma _* ) + (\xi  - \theta_* )P(\mu _2 )v_1 (\gamma _* ) =  s_* \xi u_1 (\gamma _* ).}
\end{equation}

Multiplying (\ref{n1}) by $\xi$, subtracting it from (\ref{n2}) yields $(\theta_*-\xi) \left(P(\mu_1)-P(\mu_2)\right)v_1(\gamma_*)=0.$ This is in contradiction because $v_1(\gamma_*)$ is a nonzero vector.

\textit{Case 2.}\ Suppose $\xi=\theta_*.$ Note that $\theta_*\ne 0$. In this case from (\ref{fs}) we have
\begin{equation}
u_1 (\gamma _* )^* P(\mu _1 ) + \left| \theta_*  \right|^2 u_1 (\gamma _* )^* \left( {P(\mu _1 ) - P(\mu _2 )} \right) = s_* v_1 (\gamma _* ),
\end{equation}\label{mm1}
and
\begin{equation}\label{mm2}
\bar \theta_* u_1 (\gamma _* )^* P(\mu _2 ) = s_* \bar \theta_* v_1 (\gamma _* ),
\end{equation}
Dividing (\ref{mm1}) by $\bar \theta_*$, subtracting it from (\ref{mm2}) leads to $\left(1+ \left| \theta_*  \right|^2\right) u_1 (\gamma _* )^* \left( {P(\mu _1 ) - P(\mu _2 )} \right) =0.$ This contradicts the fact that $u_1(\gamma_*)$ is a nonzero vector. \qquad$\Box$


The next corollary follows immediately.
\begin{corollary}
If $\gamma_*>0$ and $P[\mu _1 ,\mu _2 ]$ is a nonsingular matrix, then rank$(\hat V(\gamma _* )) = 2$.
\end{corollary}

\section{Computation bounds for $\mathcal{D} (P,\mu _1 ,\mu _2 ) $ and construction a perturbation of $P(\lambda)$}\label{boundsandperturbation}

In this section, at first a lower bound  of $ \mathcal{D} (P,\mu _1 ,\mu _2 )$ is computed. Then an upper bound of $\mathcal{D}(P,\mu _1 ,\mu _2 )$ will be obtained by constructing an associated perturbation of $P(\lambda)$.

\begin{lemma}\label{lowerbound}
Let $\mu_1$ and $\mu_2$ be two eigenvalues of the perturbation matrix polynomial $ Q(\lambda ) = P(\lambda ) + \Delta (\lambda
)\in{\mathcal{B}}(P,\varepsilon ,w)$. Then for any $\gamma \neq 0$
\begin{equation}\label{karan}
\varepsilon  \ge \frac{\left\| {F[\Delta
(\mu _1 ,\mu _2 );\gamma ]} \right\|}{{\left\| {\left[ {\begin{array}{*{20}c}
   {w(\left| {\mu _1 } \right|)} & 0  \\
   {\gamma \left| {w[\mu _1 ,\mu _2 ]} \right|} & {w(\left| {\mu _2 } \right|)}  \\
\end{array}} \right]} \right\|}} \ge \frac{{s_{2n - 1} (F[P(\mu _1 ,\mu _2 );\gamma ])}}{{\left\| {F[\left|w(\mu _1 ,\mu _2 )\right|;\gamma ])} \right\|}}.
\end{equation}
\end{lemma}
\textbf{Proof.} At first we have
\begin{equation*}
\left\| {\Delta [\mu _1 ,\mu _2 ]} \right\| \le \sum\limits_{j = 1}^n {\left\| {\Delta _j } \right\|} \left\| {\frac{{\left( {\mu _1^j  - \mu _2^j } \right)}}{{\mu _1  - \mu _2 }}} \right\| \le
\sum\limits_{j = 1}^n {\varepsilon w_j } \left\| {\frac{{\left( {\mu _1^j  - \mu _2^j } \right)}}{{\mu _1  - \mu _2 }}} \right\| = \varepsilon \left| {w[\mu _1 ,\mu _2 ]} \right|,
\end{equation*}
and
\begin{equation*}
\left\| {\Delta (\mu _i )} \right\| = \left\| {\sum\limits_{j = 1}^n {\Delta _j \mu _i^j } } \right\| \le \sum\limits_{j = 1}^n {\left\| {\Delta _j } \right\|\left| {\mu _i^j } \right| \le \sum\limits_{j = 1}^n {\varepsilon w_j \left| {\mu _i^j } \right| = } } \varepsilon w(\left| {\mu _i } \right|);\qquad i=1,2.
\end{equation*}
We can assume a unit vector $\left[
{\begin{array}{*{20}c}
   x  \\
   y  \\
\end{array}} \right] \in \mathbb{C}^{2n} (x,y \in \mathbb{C}^n )$ such that for any $\gamma \neq 0$, 
\begin{eqnarray*}
\left\| {F[\Delta (\mu _1 ,\mu _2 );\gamma ]} \right\|^2  &=& \left\| {\left[ {\begin{array}{*{20}c}
   {\Delta (\mu _1 )} & 0  \\
   {\gamma \Delta [\mu _1 ,\mu _2 ]} & {\Delta (\mu _2 )}  \\
\end{array}} \right]} \right\|^2 \\& =& \left\| {\left[ {\begin{array}{*{20}c}
   {\Delta (\mu _1 )} & 0  \\
   {\gamma \Delta [\mu _1 ,\mu _2 ]} & {\Delta (\mu _2 )}  \\
\end{array}} \right]\left[ {\begin{array}{*{20}c}
   x  \\
   y  \\
\end{array}} \right]} \right\|^2 \\& \le& \left\| {\left[ {\begin{array}{*{20}c}
   {\varepsilon w(\left| {\mu _1 } \right|)\left\| x \right\|}  \\
   {\left| \gamma  \right|\varepsilon \left| {w[\mu _1 ,\mu _2 ]} \right|\left\| x \right\| + \varepsilon w(\left| {\mu _2 } \right|)\left\| y \right\|}  \\
\end{array}} \right]} \right\|^2 \\& \le& \varepsilon ^2 \left\| {\left[ {\begin{array}{*{20}c}
   {w(\left| {\mu _1 } \right|)} & 0  \\
   {\left| \gamma  \right|\left| {w[\mu _1 ,\mu _2 ]} \right|} & {w(\left| {\mu _2 } \right|)}  \\
\end{array}} \right]} \right\|^2 \\& =& \varepsilon ^2 \left\| {F[\left| {w(\mu _1 ,\mu _2 )} \right|;\gamma ])} \right\|.
\end{eqnarray*}

Lemma \ref{lem2} completes this proof .\qquad $\square$

Considering Definition \ref{dis} and  Lemma \ref{lowerbound}, a lower bound for $\mathcal{D}(P,\mu _1 ,\mu _2 )$ can be obtained by minimizing the both sides of (\ref{karan}) as follows
\begin{equation}\label{Lbound}
\mathcal{D}(P,\mu _1 ,\mu _2 ) \ge \frac{{s_{2n - 1} (F[P(\mu _1 ,\mu _2 );\gamma ])}}{{\left\| {F[\left|w(\mu _1 ,\mu _2 )\right|;\gamma ])} \right\|}}.
\end{equation}
Let us now construct a perturbation of $P(\lambda)$. First assume that $\gamma_*>0$ and $P[\mu _1 ,\mu _2 ]$ is a nonsingular
matrix. Therefore, Lemma \ref{lemma7} implies that rank$\left(\hat{V}(\gamma_*)\right)=2$. In this case, a matrix polynomial
$\Delta _{\gamma_*}  (\lambda )$  is constructed such that $\mu_1$ and
$\mu_2$ are the eigenvalues of the perturbation matrix polynomial
$Q_{{\gamma_*}} (\lambda ) = P(\lambda ) + \Delta _{{\gamma_*}} (\lambda )$.

For this, define the matrix
\begin{equation}\label{deltagamma}
\Delta _{{\gamma_*}}   =  - s_*\hat{U}({\gamma_*} )\left[
{\begin{array}{*{20}c}
   {\frac{2}{{1 + \alpha _1 }}} & 0  \\
   0 & {\frac{2}{{1 + \alpha _2 }}}  \\
\end{array}} \right]\hat{V}({\gamma_*} )^\dag,
\end{equation}
where $\hat{V}({\gamma_*})^\dag$ is the
\emph{Moore-Penrose pseudoinverse} of $\hat{V}({\gamma_*} )$ and
\begin{equation*}
\alpha _1  = \frac{1}{{w(\left| {\mu _2 } \right|)}}\sum\limits_{j = 0}^m {\left( {(\frac{{\bar \mu _2 }}{{\left| {\mu _2 } \right|}})^j \mu _1^j \omega _j } \right)},\qquad {  \mbox{and}}\qquad \alpha _2  = \frac{1}{{w(\left| {\mu _1 } \right|)}}\sum\limits_{j = 0}^m {\left( {(\frac{{\bar \mu _1 }}{{\left| {\mu _1 } \right|}})^j \mu _2^j \omega _j } \right)}.
\end{equation*}
Finally the $n \times n$ matrix polynomial
$\Delta _{{\gamma_*}}  (\lambda ) = \sum\limits_{j = 0}^m {\Delta
_{{\gamma_*} ,j} \lambda ^j },$ is defined as follows
\begin{equation}\label{deltaj}
\Delta _{{\gamma_*} ,j}  = \frac{1}{2}\left(\frac{1}{{w(\left| {\mu _1
} \right|)}}(\frac{{\bar \mu _1 }}{{\left| {\mu _1 } \right|}})^j
+ \frac{1}{{w(\left| {\mu _2 } \right|)}}(\frac{{\bar \mu _2
}}{{\left| {\mu _2 } \right|}})^j \right)\omega _j \Delta _{\gamma_*},
\end{equation}
such that satisfies $\Delta_{\gamma_*} (\mu _i ) = \left(\frac{{1+ \alpha _i }}{2}\right)\Delta _{{\gamma_*}},(i=1,2).$ Keeping in mind that $u_1 ({\gamma_*} ), v_1 ({\gamma_*} ), \hat{u} ({\gamma_*} )$ and $\hat{v}({\gamma_*})$ were defined in Definition \ref{UV}, and satisfied Lemma \ref{lemma7}. For matrix polynomial
\begin{equation}\label{qgama}
Q_{\gamma_*}  (\lambda ) = P(\lambda ) + \Delta _{\gamma_*}  (\lambda ) =
\sum\limits_{j = 0}^m {(A_j  + \Delta _{{\gamma_*} ,j} )\lambda ^j
},
\end{equation}
we can obtain the following relations
\begin{eqnarray*}
Q_{\gamma _* } (\mu _1 )v_1 (\gamma _* ) &=& P(\mu _1 )v_1 (\gamma _* ) + \Delta _{\gamma _* } (\mu _1 )v_1 (\gamma _* ) \\&=& s_* u_1 (\gamma _* ) + \left( {\frac{{1 + \alpha _1 }}{2}} \right)\Delta _{\gamma _* } v_1 (\gamma _* ) \\&=& s_* u_1 (\gamma _* ) + \left( {\frac{{1 + \alpha _1 }}{2}} \right)\left( {\frac{2}{{1 + \alpha _1 }}} \right)\left( - s_* u_1 (\gamma _* )\right)\\& =& 0,
\end{eqnarray*}
and
\begin{eqnarray*}
Q_{\gamma _* } (\mu _2 )\hat v(\gamma _* )& =& P(\mu _2 )\hat v(\gamma _* ) + \Delta _{\gamma _* } (\mu _2 )\hat v(\gamma _* ) \\&=& s_* \hat u(\gamma _* ) + \left( {\frac{{1 + \alpha _2 }}{2}} \right)\Delta _{\gamma _* } \hat v(\gamma _* )\\& =& s_* \hat u(\gamma _* ) + \left( {\frac{{1 + \alpha _2 }}{2}} \right)\left( {\frac{2}{{1 + \alpha _2 }}} \right)\left( - s_* \hat u(\gamma _* )\right) \\&=& 0.
\end{eqnarray*}
Consequently $\mu_1$ and $\mu_2$ are two eigenvalues of $
Q_{\gamma_*} (\lambda )$ corresponding to $v_1 ({\gamma_*} ) $ and
$\hat{v} ({\gamma_*} ) $ as two eigenvectors, respectively. On the other hand, it follows from (\ref{deltaj}) that
\begin{eqnarray*}
\left\| {\Delta _{{\gamma_*} ,j} } \right\| \le \frac{{\omega
_j }}{2}(\frac{1}{{w(\left| {\mu _1 } \right|)}} +
\frac{1}{{w(\left| {\mu _2 } \right|)}}) \left\| {\Delta _{\gamma_*}  } \right\|, \qquad j =
0,1,\ldots,m.
\end{eqnarray*}
Consequently, an upper bound of $ \mathcal{D}(P,\mu _1 ,\mu _2 )$ is obtained the by following relation for any ${\gamma_*} > 0$
\begin{equation}\label{Ubound}
\mathcal{D} (P,\mu _1 ,\mu _2 ) \le \frac{{1
}}{2}\left(\frac{1}{{w(\left| {\mu _1 } \right|)}} +
\frac{1}{{w(\left| {\mu _2 } \right|)}}\right)\left\| {\Delta _{\gamma_*}  }
\right\|.
\end{equation}
It will be convenient to represent the lower bound provided in (\ref{Lbound}) by $\beta _{low} (P,\mu _1 ,\mu _2 ,{\gamma})$ and the upper bound provided in (\ref{Ubound}) by $\beta _{up} (P,\mu _1 ,\mu _2 ,{\gamma})$, i.e.,
\begin{equation}\label{low}
\beta _{low} (P,\mu _1 ,\mu _2 ,{\gamma}) =\frac{{s_{2n - 1} (F[P(\mu _1 ,\mu _2 );\gamma ])}}{{\left\| {F[\left|w(\mu _1 ,\mu _2 )\right|;\gamma ])} \right\|}},
\end{equation}
and
\begin{equation}\label{up}
\beta _{up} (P,\mu _1 ,\mu _2 ,{\gamma_*}) = \frac{{1
}}{2}\left(\frac{1}{{w(\left| {\mu _1 } \right|)}} +
\frac{1}{{w(\left| {\mu _2 } \right|)}}\right)\left\| {\Delta _{\gamma_*}  }
\right\|.
\end{equation}
The results obtained so far from the beginning of this section are summarized in the next theorem.
\begin{theorem}\label{thm11}
Let $P(\lambda)$ be the matrix polynomial as in  {(\ref{plambda})} and let
$\mu_1$ and $\mu_2$ be two given distinct complex numbers. Then for any
$\gamma>0$,
\begin{equation*}
\beta _{low} (P,\mu _1 ,\mu _2 ,\gamma) \le \mathcal{D}(P,\mu _1
,\mu _2 ),
\end{equation*}
where $\beta _{low} (P,\mu _1 ,\mu _2 ,\gamma)$ is introduced in  {(\ref{low})}. In addition, if $\gamma_* > 0$,
then the matrix polynomial $ Q_{\gamma_*}
(\lambda )$ in (\ref{qgama}) has  $\mu_1$ and $\mu_2$ as two its eigenvalues corresponding to $v_1 (\gamma_* ) $ and $\hat{v} (\gamma_* ) $ as two its eigenvectors, respectively. Furthermore, $Q_{\gamma_*}
(\lambda ) \in \partial {\mathcal{B}}(P,\beta _{up} (P,\mu _1 ,\mu _2
,\gamma_* ),w)$ and $\mathcal{D}(P,\mu _1 ,\mu _2 ) \le \beta _{up}
(P,\mu _1 ,\mu _2 ,\gamma_* ),$
where $\beta _{up} (P,\mu _1 ,\mu _2 ,\gamma_*)$ is introduced in (\ref{up}).
\end{theorem}
It should be pointed out that the bounds obtained are not necessarily optimal, however, it is assured that $\mathcal{D}(P,\mu _1 ,\mu _2 )$ belongs to $ \left[ {\beta _{low}(P,\mu _1 ,\mu _2 ,\gamma _* ),\beta _{up} (P,\mu _1 ,\mu _2,\gamma _* )} \right] $. Anyhow, the following remark can be used to obtain some close bounds.  
\begin{remark}\label{fminn}
It is important to note that, Theorem \ref{thm11} holds for any $\gamma_0>0$ that assures the matrix $\hat{V}(\gamma_0)=$ is full (column) rank. Therefore, it can be an obvious expectation to find a value of $\gamma>0$ that obtains the closest upper and lower bounds. For doing this, we can define the following nonnegative function
\[f(\gamma)=\beta _{up} (P,\mu _1 ,\mu _2 ,{\gamma})-\beta _{low} (P,,\mu _1 ,\mu _2 ,{\gamma}),\]
and try to minimize this function by implementation of unconstrained optimization methods (for example, see \cite{nocedal}). On the other hand, best lower bound and finest upper bound can be obtained by maximizing and minimizing $\beta _{low} (P,\mu _1 ,\mu _2 ,{\gamma})$ and $\beta _{up} (P,\mu _1 ,\mu _2 ,{\gamma})$, respectively. It is clear that values of $\gamma$ which yield the smallest upper bound and the biggest lower bound may be different.
\end{remark}
Now suppose that the singular value $s_{2n - 1} (F[P(\mu _1 ,\mu
_2 );\gamma ])$ attains its maximum value at $\gamma=0,$ i.e., $\gamma_*=0.$ Next we compute an upper bound for $\mathcal{D}(P,\mu _1 ,\mu _2 ) \le \beta _{up}
(P,\mu _1 ,\mu _2 ,\gamma_* ),$ constructing associated perturbations of $P(\lambda)$.

Let $u_i ,v_i \in \mathbb{C}^n,~(i=1,2)$ be a pair of left and right
singular vectors of $P(\mu_i)$ corresponding to
$\sigma_i=s_n(P(\mu_i)),~(i=1,2)$, respectively, such that $v_1$ and $v_2$ are linearly independent. We define the matrix polynomial $\Delta_0(\lambda)$ as
\begin{equation}
\Delta _0 (\lambda ) = \Delta _0  =  - \left[
{\begin{array}{*{20}c}
   {u_1 } & {u_2 }  \\
\end{array}} \right]\left[ {\begin{array}{*{20}c}
   {\sigma_1 } & 0  \\
   0 & {\sigma_2}  \\
\end{array}} \right]\left[ {\begin{array}{*{20}c}
   {v_1 } & {v_2 }  \\
\end{array}} \right]^\dag,
\end{equation}
where $\left[ {\begin{array}{*{20}c}
   {v_1 } & {v_2 }  \\
\end{array}} \right]^\dag$ is the Moore-Penrose pseudoinverse of
$
\left[ {\begin{array}{*{20}c}
   {v_1 } & {v_2 }  \\
\end{array}} \right]
$.
Then,
\begin{equation}\label{q0}
Q_0 (\lambda ) = P(\lambda ) + \Delta _0 (\lambda ) = A_m \lambda
^m  + A_{m - 1} \lambda ^{m - 1}  + ... + A_1 \lambda  + \left(A_0
+ \Delta _0 \right) ,
\end{equation}
lies on $
\partial {\mathcal{B}}(P,\frac{{\left\| {\Delta _0 } \right\|}}{{\omega _0 }},\varepsilon )
$ and satisfies
\[
Q_0 (\mu _i )v_i  = P(\mu _i )v_i  + \Delta _0 (\mu _i )v_i  =
\sigma_i u_i  - \sigma_i u_i  = 0;\qquad i = 1,2.
\]
Hence $\mu_1$ and $\mu_2$ are two eigenvalues of the matrix
polynomial $Q_0 (\lambda)$ with corresponding eigenvectors $v_1$
and $v_2$, respectively.
\begin{theorem}
Let $\gamma_*=0$, and let $u_i ,v_i  \in \mathbb{C}^n,~(i=1,2)$ be a pair of left and right singular vectors of
$P(\mu_i)$ corresponding to $\sigma_i=s_n(P(\mu_i)),~(i=1,2)$, respectively. If $v_1$ and $v_2$ are linearly independent, then the matrix polynomial $Q_0 (\lambda)$ in (\ref{q0})
lies on $\partial \mathcal{B}(P,\frac{{\left\| {\Delta _0 } \right\|}}{{\omega _0 }},\varepsilon )$ and has $\mu_1$ and $\mu_2$ as its eigenvalues associated with $v_1, v_2$ as two its eigenvectors, respectively.
\end{theorem}
Two special cases of the matter of our discussion, are considered in two following remarks.
\begin{remark}\label{rimark}
Suppose that we have  $\omega _0>0$ and $ \omega _1  = \omega _2
= \cdots = \omega _m  = 0 $ as a given set of nonnegative
weights. Then $\alpha_1=\alpha_2=1, \omega (\left| {\mu _1 }
\right|) = \omega (\left| {\mu _2 } \right|) = \omega _0$ and thus $
\hat \omega [\mu _1 ,\mu _2 ] = 0$. Consequently, in this case the lower and
upper bounds of $\mathcal{D}(P,\mu _1 ,\mu _2 )$ given by (\ref{low}) and (\ref{up}) respectively, are equal to $
\frac{s_*}{{\omega _0 }}$ and thus $\mathcal{D}(P,\mu _1 ,\mu _2 )$ is determined exactly.
\end{remark}
\begin{remark}
Assume that $P(\lambda ) = I\lambda  - A,$ where $A \in \mathbb{C}^{n\times n} $ and $w=\{\omega_0, \omega_1\}=\{1, 0\}$. Then $F[P(\mu _1 ,\mu _2 );\gamma ]$ will be coincides with the matrix in the results \cite{gracia}, i.e., 
\begin{equation*}
 F[P(\mu _1 ,\mu _2 );\gamma ] = \left[
{\begin{array}{*{20}c}
   {P(\mu_1 )} & 0  \\
   {\gamma P[\mu _1 ,\mu _2 ]} & {P(\mu_2 )}  \\
\end{array}} \right] = \left[ {\begin{array}{*{20}c}
   {I\mu_1  - A} & 0  \\
   {\gamma I} & {I\mu_2  - A}  \\
\end{array}} \right],
\end{equation*}
Also, one can find that $w(\lambda ) = 1$ and $w'(\lambda ) =
0$, which yield $\alpha_1=\alpha_2=1.$
Thus if $\gamma_* >0$, then the matrix $\Delta _{\gamma_{*}}$ in (\ref{deltagamma}) becomes $\Delta _{\gamma_{*}}   =  - s_* \hat{U}(\gamma_*) \hat{V}(\gamma_*)^\dag$. By Corollary \ref{natijeh}, the two matrices $\hat{U}(\gamma _* )$ and $\hat{V}(\gamma _* )$ have the same nonzero singular values. Therefore, there exists a  unitary matrix $W \in \mathbb{C}^{n \times n}$ such that $\hat{U}(\gamma_*)=W \hat{V}(\gamma_*)$. In addition, the fact that $\hat{V}(\gamma _* )$ is a full rank matrix concludes $\hat{V}(\gamma_* )^{\dag} \hat{V}(\gamma_* )=I_k$. Therefore, it follows
\begin{equation*}
\mathcal{D}(P,\mu_1,\mu_2)={\left\| \Delta_{\gamma_*}  \right\|_2} = {\left\| { - {s_*}\hat{U}({\gamma _*}){\hat{V} }({\gamma _*})^\dag} \right\|_2} =s_*{\left\| {W\hat{V}({\gamma _*}){\hat{V} }({\gamma _*})^\dag} \right\|_2} = s_*.
\end{equation*}

Furthermore, the perturbation matrix
polynomial $Q_{\gamma _* }(\lambda)$ in (\ref{qgama}) turns into
\[ Q_{\gamma _* } (\lambda ) = I\lambda  - \left( {A + s_* \hat
U(\gamma _* )\hat V(\gamma _* )} \right).
\]
Consequently, the results obtained in this article can be construed as a generalization of the results obtained in \cite{lipert, gracia} for the case of matrix polynomials.
\end{remark}
\section{Numerical experiments }
Review the topic of this paper and let us concentrate on the subject of finding a matrix polynomial that has two eigenvalues located at desired positions. This viewpoint can by useful in some problems such as reconstructing a matrix polynomial from prescribed spectral data which can be assumed as inverse eigenvalue problem for the case of matrix polynomials. Assume now, we are asked to find a matrix polynomial having to given scalars $\mu_1, \mu_2 \in \mathbb{C}$ as its eigenvalues. For doing this, one can consider an arbitrary matrix polynomial, namely, $P(\lambda)$ in the craved size. Next, by following procedure the described in Section 4, the desired matrix polynomial (which $\mu_1$ and $\mu_2$ are some of its eigenvalues) is computable. All computations were performed in Matlab with 16 significant figures, however, for simplicity all numerical results are shown with 4 decimal places.
\begin{example}
Let two scalars $\mu_1=1$ and $\mu_2=2+i$ are given and let we are asked to find a $2\times 2$ matrix polynomial such that $\mu_1$ and $\mu_2$ are some of its eigenvalues. To do this, consider
\[
P(\lambda ) = \left[ {\begin{array}{*{20}c}
   -5 & 10  \\
   -4 & 5  \\
\end{array}} \right]\lambda ^2  + \left[ {\begin{array}{*{20}c}
   { - 2} & { - 4}  \\
   1 & { - 1}  \\
\end{array}} \right]\lambda  + \left[ {\begin{array}{*{20}c}
   -1 & -3  \\
   0 & 0  \\
\end{array}} \right],
\]
which its coefficient matrices are randomly generated by MATLAB. Let the set of weights $ w = \{3.1623, 4.4966, 12.8310\}$ be the norms of the coefficient matrices. Employing the MATLAB function \texttt{fminsearch} we find that
$\gamma_{*}=1.8914$ and 
$s_{*}=s_{3}(F[P(1,2+i);\gamma_*])=4.1132.$
The graph of the $s_{3}(F[P(1,2+i);\gamma])$ for $\gamma\in [0,10]$ is plotted in Fig 1 and $(\gamma_* ,s_{*})$ is marked with "o". Now, applying the procedures described in Section \ref{boundsandperturbation}, we can compute the matrix polynomial $Q_{1.8914} (\lambda ) = P(\lambda)+\Delta _{1.8914}(\lambda)$ as a perturbation of $P(\lambda)$ that lies on $
\partial {\mathcal{B}}(P,\beta _{up} (P,1,2+i ,1.8914),w)$ and has $\mu_1=1$ and $\mu_2=2+i$ as two eigenvalues. Where
\begin{eqnarray*}
\Delta _{1.8914} (\lambda )&  = & \left[ {\begin{array}{*{20}c}
   {{ {0}}{ {.4834  -  0}}{ {.6940i}}} & {{ { - 1}}{ {.2959  +  0}}{ {.5336i}}}  \\
   {{ {1}}{ {.8038  +  0}}{ {.2529i}}} & {{ { - 2}}{ {.4162  +  0}}{ {.0769i}}}  \\
\end{array}} \right]\lambda ^2 \\& + & \left[ {\begin{array}{*{20}c}
   {{ {0}}{ {.1999  -  0}}{ {.2403i}}} & {{ { - 0}}{ {.4941  +  0}}{ {.1557i}}}  \\
   {{ {0}}{ {.6563  +  0}}{ {.1500i}}} & {{ { - 0}}{ {.8922  -  0}}{ {.0478i}}}  \\
\end{array}} \right]\lambda \\&  + & \left[ {\begin{array}{*{20}c}
   {{ {0}}{ {.1588  -  0}}{ {.1576i}}} & {{ { - 0}}{ {.3627  +  0}}{ {.0772i}}}  \\
   {{ {0}}{ {.4575  +  0}}{ {.1517i}}} & {{ { - 0}}{ {.6326  -  0}}{ {.0949i}}}  \\
\end{array}} \right].
\end{eqnarray*}

In addition, we can calculate lower and upper bounds for $\mathcal{D}(P, 1,2+i)$ as follows

$\beta _{{ {low}}} ({ {P}}, 1,2+i,1.8914) = { {0}}.{ {0376}}\,\qquad {  \mbox{and}}\qquad \beta _{{ {up}}} ({ {P}}, 1,2+i, 1.8914) = { {0}}.{ {2847}}.$
\begin{figure}
 \centering
 \includegraphics[width=0.53\linewidth]{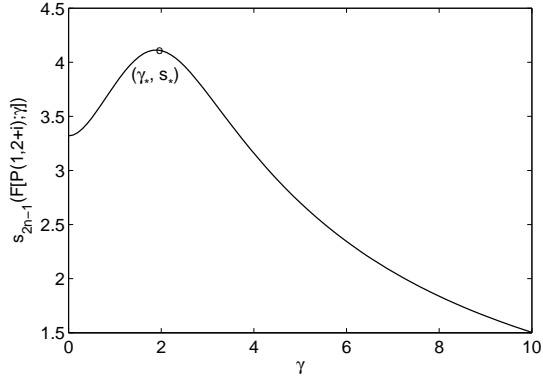}
 \caption{The graph of the $s_{3}(F[P(1,2+i);\gamma])$.}
 \label{fig:ex1}
 \end{figure}
The graphs of bounds $\beta _{up} (P,1,2+i,\gamma )$ and $\beta _{low} (P,1,2+i ,\gamma )$
are plotted in Fig 2, for $\gamma\in [0, 10]$ and the bounds $\beta _{up}(P,1,2+i, 1.8914)$ and $\beta _{low} (P,1,2+i, 1.8914)$ are
marked with "o".
\begin{figure}
   \centering
   \includegraphics[width=0.53\linewidth]{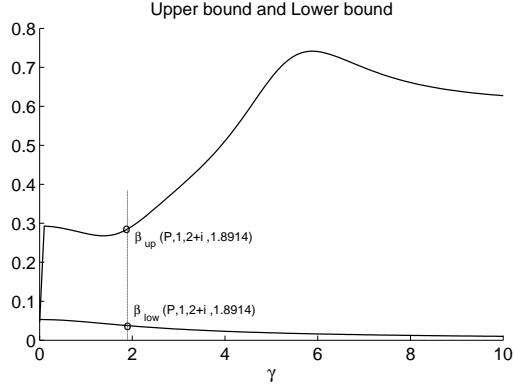}
   \caption{The graphs of the bounds $\beta _{up} (P,1,2+i, \gamma )$ and $\beta _{low} (P,1,2+i ,\gamma )$.}
   \label{fig:karanhanew}
   \end{figure}
Furthermore, next results verify Lemma \ref{lemma6} and Corollary \ref{natijeh}, respectively,
\begin{equation*}
\left| {u_2^* (\gamma _* )P[\mu _1 ,\mu _2 ]v_1 (\gamma _*
)} \right| = { {1}}{ {.5612}} \times { {10}}^{ - 5},
\end{equation*}
\begin{equation*}
\left\| {\hat{U} (\gamma _* )^*\hat{U}(\gamma _* ) - \hat{V} (\gamma _* )^*\hat{V}(\gamma
_* )} \right\| = { {2}}{ {.9886}} \times { {10}}^{ - 5}.
\end{equation*}
\end{example}
We also consider an example for the case of $\gamma_*=0$.
\begin{example}
Let
\[
P(\lambda ) = \left[ {\begin{array}{*{20}c}
   3 & 0 & 1  \\
   8 & { - 1} & 0  \\
   4 & 2 & 3  \\
\end{array}} \right]\lambda ^2  + \left[ {\begin{array}{*{20}c}
   6 & { - 4} & 0  \\
   1 & { - 5} & 5  \\
   1 & { - 1} & {10}  \\
\end{array}} \right]\lambda  + \left[ {\begin{array}{*{20}c}
   9 & 7 & 6  \\
   2 & 7 & { - 4}  \\
   { - 2} & 6 & 5  \\
\end{array}} \right],
\]
and let $w = \{ 1,1,1\}$. Now consider two
complex numbers $\mu_1=5$ and $\mu_2=-1$. It is easy to find that $s_{5}(F[P(5,-1);\gamma])$ attains its maximum at $\gamma_* = 0$ and $s_*=4.0378$.
\begin{figure}
\centering
\includegraphics[width=0.53\linewidth]{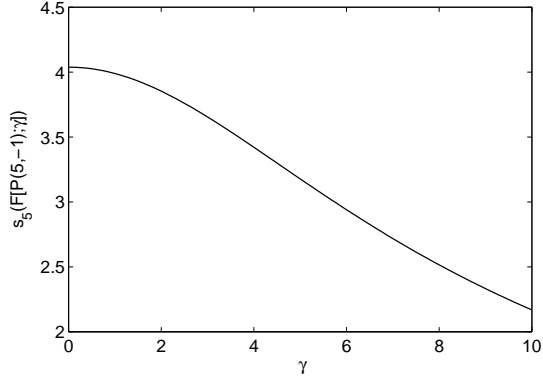}
\caption{The graph of the $s_{5}(F[P(5,-1);\gamma])$.}
\label{fig:gammasefr}
\end{figure} 
The graph of the $s_{5}(F[P(5,-1);\gamma])$ for $\gamma\in [0,
10]$ is plotted in Fig 3.
According to discussion for the case
$\gamma_*=0$, we obtain the matrix polynomial $Q_0 (\lambda
)=P(\lambda )+\Delta _0$ belonging to $
\partial \mathcal{B}(P,{ {4}}{{.1545,w)}}$ having $\mu_1=5$ and
$\mu_2=-1$ as its two eigenvalues. Where
\[
\Delta _0 (\lambda ) = \Delta _0  = \left[ {\begin{array}{*{20}c}
   {{ { - 0}}{ {.6257}}} & {{ {0}}{ {.8167}}} & {{ { - 0}}{ {.3709}}}  \\
   {{ { - 1}}{ {.5026}}} & {{ {0}}{ {.0959}}} & {{ {0}}{ {.1659}}}  \\
   {{ {3}}{ {.6390}}} & {{ { - 1}}{ {.0783}}} & {{ {0}}{ {.0774}}}  \\
\end{array}} \right].
\]
\end{example}
\section{Conclusions}
In this paper, for a matrix polynomial $P(\lambda)$ and two given
distinct complex numbers $\mu_1$ and $\mu_2$, a spectral norm
distance from $P(\lambda)$ to the set of matrix polynomials that
have $\mu_1$ and $\mu_2$ as two eigenvalues, was introduced. The
upper and lower bounds for this distance were computed and
associated perturbation of $P(\lambda)$ was constructed. The cases
of $\gamma_*>0$ and $\gamma_*=0$ were studied in detail
separately. Finally, it was pointed out that the bounds obtained
are not necessarily optimal, however, it is assured that $\mathcal{D}(P,\mu _1 ,\mu _2 )$ belongs to $ \left[ {\beta _{low}
(P,\mu _1 ,\mu _2 ,\gamma _* ),\beta _{up} (P,\mu _1 ,\mu _2
,\gamma _* )} \right] $. The conditions to obtain the optimal bounds is the subject of our future research.

\end{document}